# OPTIMAL DESIGNS FOR A CLASS OF NONLINEAR REGRESSION MODELS

By Holger Dette[1], Viatcheslav B. Melas and Andrey Pepelyshev

*Ruhr-Universität Bochum, St. Petersburg State University and St. Petersburg State University*

For a broad class of nonlinear regression models we investigate the local $E$- and $c$-optimal design problem. It is demonstrated that in many cases the optimal designs with respect to these optimality criteria are supported at the Chebyshev points, which are the local extrema of the equi-oscillating best approximation of the function $f_0 \equiv 0$ by a normalized linear combination of the regression functions in the corresponding linearized model. The class of models includes rational, logistic and exponential models and for the rational regression models the $E$- and $c$-optimal design problem is solved explicitly in many cases.

**1. Introduction.** Nonlinear regression models are widely used to describe the dependencies between a response and an explanatory variable [see, e.g., Seber and Wild (1989), Ratkowsky (1983) or Ratkowsky (1990)]. An appropriate choice of the experimental conditions can improve the quality of statistical inference substantially and, therefore, many authors have discussed the problem of designing experiments for nonlinear regression models. We refer to Chernoff (1953) and Melas (1978) for early references and Ford, Torsney and Wu (1992), He, Studden and Sun (1996) and Dette, Haines and Imhof (1999) for more recent references on local optimal designs. Because local optimal designs depend on an initial guess for the unknown parameters, several authors have proposed alternative design strategies. Bayesian or robust optimal designs have been discussed by Pronzato and Walter (1985) and Chaloner and Larntz (1989), among many others [see Chaloner and Verdinelli

Received December 2002; revised November 2003.
[1]Supported by the Deutsche Forschungsgemeinschaft (SFB 475: Komplexitätsreduktion in multivariaten Datenstrukturen, Sachbeihilfe De 502/18-1, 436 RUS 113/712/0-1).

(1995) and the references therein]. Other authors propose sequential methods, which update the information about the unknown parameter sequentially [see, e.g., Ford and Silvey (1980) and Wu (1985)]. Most of the literature concentrates on $D$-optimal designs (independent of the particular approach), which maximize the determinant of the Fisher information matrix for the parameters in the model, but much less attention has been paid to $E$-optimal designs in nonlinear regression models, which maximize the minimum eigenvalue of the Fisher information matrix [see Dette and Haines (1994) or Dette and Wong (1999), who report some results for models with two parameters].

Because local optimal designs are the basis for all advanced design strategies, it is the purpose of the present paper to study local $E$-optimal designs for a class of nonlinear regression models which can be represented in the form

$$(1.1) \quad Y = \sum_{i=1}^{s} a_i h_i(t) + \sum_{i=1}^{k} a_{s+i} \varphi(t, b_i) + \varepsilon.$$

Here $\varphi$ is a known function, the explanatory variable $t$ varies in an interval $I \subset \mathbb{R}$, $\varepsilon$ denotes a random error with mean zero and constant variance and $a_1, \ldots, a_{s+k}, b_1, \ldots, b_k \in \mathbb{R}$ denote the unknown parameters of the model. The consideration of this type of model was motivated by recent work of Imhof and Studden (2001), who considered a class of rational models of the form

$$(1.2) \quad Y = \sum_{i=1}^{s} a_i t^{i-1} + \sum_{i=1}^{k} \frac{a_{s+i}}{t - b_i} + \varepsilon,$$

where $t \in I, b_i \neq b_j$ $(i \neq j)$ and the parameters $b_i \notin I$ are assumed to be known for all $i = 1, \ldots, k$. Note that model (1.2) is in fact linear, because Imhof and Studden (2001) assumed the $b_i$ to be known. These models are very popular because they have appealing approximation properties [see Petrushev and Popov (1987) for some theoretical properties and Dudzinski and Mykytowycz (1961) and Ratkowsky (1983), page 120, for an application of this model]. In this paper [in contrast to the work of Imhof and Studden (2001)] the parameters $b_1, \ldots, b_k$ in the model (1.1) are not assumed to be known, but also have to be estimated from the data. Moreover, the model (1.1) considered here includes numerous other regression functions. For example, in environmental and ecological statistics exponential models of the form $a_1 e^{b_1 t} + a_2 e^{b_2 t}$ are frequently used in toxicokinetic experiments [see, e.g., Becka and Urfer (1996) or Becka, Bolt and Urfer (1993)] and this corresponds to the choice $\varphi(t, x) = e^{tx}$ in (1.1). Another popular class of logarithmic models is obtained from equation (1.1) by the choice $\varphi(t, x) = \log(t - x)$.

Imhof and Studden (2001) studied $E$-optimal designs for the model (1.2) with $s = 1$ under the assumption that the nonlinear parameters $b_1, \ldots, b_k$ are known by the experimenter and do not have to be estimated from the data. In particular, they proved that the support of the $E$-optimal design for estimating a subset of the parameters $a_1, \ldots, a_{\ell+1}$ is given by the Chebyshev points corresponding to the functions $1, \frac{1}{t-b_1}, \ldots, \frac{1}{t-b_k}$ in the model (1.2). These points are the extremal points of the function $1 + \sum_{i=1}^{k} \frac{a_i^*}{x-b_i} = p^*(x)$, in the interval $I$, which has the smallest deviation from zero, that is,

$$(1.3) \qquad \sup_{x \in I} |p^*(x)| = \min_{a_2, \ldots, a_{k+1}} \sup_{x \in I} \left| 1 + \sum_{i=1}^{k} \frac{a_i}{x-b_i} \right|.$$

The universality of this solution is due to the fact that any subsystem of the regression functions in the model (1.2), which is obtained by deleting one of the basis functions, forms a weak Chebyshev system on the interval $I$ [see Karlin and Studden (1966) or the discussion in Section 2]. However, in the case where parameters $b_1, \ldots, b_k$ are unknown and also have to be estimated from the data, the local optimal design problem for the model (1.2) is equivalent to an optimal design problem in the linear regression model

$$(1.4) \qquad Y = \sum_{i=1}^{s} \beta_i t^{i-1} + \sum_{i=1}^{k} \left( \frac{\beta_{s+2i-1}}{t-b_i} + \frac{\beta_{s+2i}}{(t-b_i)^2} \right) + \varepsilon,$$

for which the corresponding regression functions do not satisfy the weak Chebyshev property mentioned above. Nevertheless, we will prove in this paper that in cases with $k \geq 2$, where the quantity $\max_{i \neq j} |b_i - b_j|$ is sufficiently small, local $E$-optimal designs and many local $c$-optimal designs for estimating linear combinations of the parameters are still supported on Chebyshev points. This fact simplifies the construction of local $E$-optimal designs substantially. Moreover, we show that this result does not depend on the specific form of the model (1.2) and (1.4) but can be established for the general model (1.1) (or its equivalent linearized model). Additionally, it can be shown numerically that in many cases the $E$-optimal design is, in fact, supported on the Chebyshev points for all admissible values of the parameters $b_1, \ldots, b_k$ ($b_i \neq b_j; i \neq j$). Our approach is based on a study of the limiting behavior of the information matrix in model (1.1) in the case where all nonlinear parameters in the model (1.1) tend to the same limit. We show that in this case the local $E$-optimal and many local optimal designs for estimating linear combinations of the coefficients $a_{s+1}, b_{s+1}, \ldots, a_{s+k}, b_{s+k}$ in the model (1.1) have the same limiting design. This indicates that $E$-optimal designs in models of type (1.1) yield precise estimates of the individual coefficients and we will illustrate this fact in several concrete examples.

It is notable that the results regarding local $E$- and $c$-optimal designs in the regression model (1.1) based on Chebyshev approximation are obtained under the simplifying assumption that $b_i = x + \delta r_i$ $(i = 1, \ldots, k)$ with $r_i \neq r_j$ and $\delta$ sufficiently small. Obviously, every vector $b = (b_1, \ldots, b_k)$ can be represented in this form, but the answer to the question if $\delta$ is sufficiently small such that our results are applicable depends on the basic function $\varphi$ used in (1.1) and the vector $b$ itself. However, the theoretical results of this paper suggest a simple procedure to obtain $E$- and $c$-optimal designs for the model (1.1). We use the designs derived under the simplifying assumption to obtain candidates for the optimal designs and check the optimality of these candidates by using equivalence theorems or alternative characterizations. Moreover, the examples of this paper and additional examples in a technical report of Dette, Melas and Pepelyshev (2002) indicate that in many cases the designs obtained under the simplifying assumption yield, in fact, the $E$- or $c$-optimal designs.

The remaining part of the paper is organized as follows. In Section 2 we introduce the basic concepts and notation, and present some preliminary results. Section 3 is devoted to an asymptotic analysis of the model (1.1), which is based on a linear transformation introduced in the Appendix. Finally, some applications to the rational model (1.2) and its equivalent linear regression model (1.4) are presented in Section 4, which extend the results of Imhof and Studden (2001) to the case where the nonlinear parameters in the model (1.2) are not known and have to be estimated from the data. Finally, all proofs and technical details are deferred to the Appendix.

**2. Preliminary results.** Consider the nonlinear regression model (1.1) and define

$$
\begin{aligned}
f(t,b) &= (f_1(t,b), \ldots, f_m(t,b))^T \\
&= (h_1(t), \ldots, h_s(t), \varphi(t,b_1), \varphi'(t,b_1), \ldots, \varphi(t,b_k), \varphi'(t,b_k))^T
\end{aligned}
\tag{2.1}
$$

as a vector of $m = s + 2k$ regression functions, where the derivatives of the function $\varphi$ are taken with respect to the second argument. It is straightforward to show that the Fisher information for the parameter $(\beta_1, \ldots, \beta_m)^T = \beta$ in the linear regression model

$$
\begin{aligned}
Y &= \beta^T f(t,b) + \varepsilon \\
&= \sum_{i=1}^{s} \beta_i h_i(t) + \sum_{i=1}^{k} (\beta_{s+2i-1} \varphi(t,b_i) + \beta_{s+2i} \varphi'(t,b_i)) + \varepsilon
\end{aligned}
\tag{2.2}
$$

is given by $f(t,b)f^T(t,b)$. Following Kiefer (1974), we call any probability measure $\xi$ with finite support on the interval $I$ an (approximate) design. The support points give the locations where observations have to be taken, while the masses correspond to the relative proportions of total observations to be taken at the particular points. For a design $\xi$ the information matrix in model (2.2) is defined by

$$(2.3) \qquad M(\xi, b) = \int_I f(t,b) f^T(t,b) \, d\xi(t),$$

and a local optimal design maximizes an appropriate function of the information matrix [see Silvey (1980) or Pukelsheim (1993)]. The dependence on the parameter $b$ is omitted whenever it is clear from the context. Among the numerous optimality criteria proposed in the literature, we consider the $E$- and $c$-optimality criteria in this paper. An $E$-optimal design $\xi_E^*$ maximizes the minimum eigenvalue $\lambda_{\min}(M(\xi, b))$ over the set of all approximate designs, while for a given vector $c \in \mathbb{R}^m$ a $c$-optimal design minimizes the expression $c^T M^-(\xi, b) c$, where the minimum is taken over the set of all designs for which the linear combination $c^T \beta$ is estimable, that is, $c \in \text{range}(M(\xi, b)) \; \forall b$. A particular case appears for the choice $c = e_i$, where $e_i \in \mathbb{R}^m$ $(i = 1, \ldots, m)$ is the $i$th unit vector. In this case we call the $c$-optimal design optimal for estimating the individual coefficient $\beta_i$.

Note that the information matrix in the nonlinear regression model (1.1) is given by $K_a^{-1} M(\xi, b) K_a^{-1}$, where the matrix $K_a \in \mathbb{R}^{m \times m}$ is defined by

$$(2.4) \qquad K_a = \text{diag}\left(\underbrace{1, \ldots, 1}_{s}, \underbrace{1, \frac{1}{a_1}, 1, \ldots, 1, \frac{1}{a_k}}_{2k}\right).$$

Consequently, a local optimal design problem in a nonlinear model (1.1) corresponds to an optimal design problem in model (2.2) for the transformed vector of parameters $K_a b$. For example, the $c$-optimal design for the model (1.1) can be obtained from the $\bar{c}$-optimal design in model (2.2), where the vector $\bar{c}$ is given by $\bar{c} = K_a c$. Similarly, the local $E$-optimal design in the nonlinear regression model (1.1) maximizes $\lambda_{\min}(K_a^{-1} M(\xi, b) K_a^{-1})$, where $M(\xi, b)$ is the information matrix in the equivalent linear regression model (2.2). For the sake of transparency we will mainly concentrate on the linearized version (2.2). The corresponding results in the nonlinear regression model (1.1) will be briefly mentioned whenever it is necessary.

A set of functions $f_1, \ldots, f_m : I \to \mathbb{R}$ is called a weak Chebyshev system (on the interval $I$) if there exists an $\varepsilon \in \{-1, 1\}$ such that

$$(2.5) \qquad \varepsilon \cdot \begin{vmatrix} f_1(x_1) & \cdots & f_1(x_m) \\ \vdots & \ddots & \vdots \\ f_m(x_1) & \cdots & f_m(x_m) \end{vmatrix} \geq 0$$

for all $x_1, \ldots, x_m \in I$ with $x_1 < x_2 < \cdots < x_m$. If the inequality in (2.5) is strict, then $\{f_1, \ldots, f_m\}$ is called a Chebyshev system. It is well known [see

Karlin and Studden (1966), Theorem II 10.2] that if $\{f_1, \ldots, f_m\}$ is a weak Chebyshev system, then there exists a unique function

$$\sum_{i=1}^{m} c_i^* f_i(t) = c^{*T} f(t) \tag{2.6}$$

with the following properties:

(2.7)
(i) $|c^{*T} f(t)| \leq 1 \ \forall t \in I$,
(ii) there exist $m$ points $s_1 < \cdots < s_m$ such that
$$c^{*T} f(s_i) = (-1)^i, \qquad i = 1, \ldots, m.$$

The function $c^{*T} f(t)$ is called a Chebyshev polynomial and the points $s_1, \ldots, s_m$ are called Chebyshev points and need not be unique. They are unique if $1 \in \text{span}\{f_1, \ldots, f_m\}, m \geq 1$, and $I$ is a bounded and closed interval, in which case $s_1 = \min_{x \in I} x$, $s_m = \max_{x \in I} x$. It is well known [see Studden (1968), Pukelsheim and Studden (1993), Heiligers (1994) or Imhof and Studden (2001), among others] that in many cases the $E$- and $c$-optimal designs in the linear regression model

$$Y = \beta^T f(t) + \varepsilon \tag{2.8}$$

are supported at the Chebyshev points. For the following discussion assume that the functions $f_1, \ldots, f_m$ generate a Chebyshev system on the interval $I$ with Chebyshev polynomial $c^{*T} f(t)$ and Chebyshev points $s_1, \ldots, s_m$, define the $m \times m$ matrix $F = (f_i(s_j))_{i,j=1}^m$ and consider a vector of weights given by

$$w = (w_1, \ldots, w_m)^T = \frac{JF^{-1}c^*}{\|c^*\|^2}, \tag{2.9}$$

where the matrix $J$ is defined by $J = \text{diag}\{(-1), 1, \ldots, (-1)^m\}$. It is then easy to see that

$$\frac{c^*}{\|c^*\|^2} = FJw = \sum_{j=1}^{m} f(s_j)(-1)^j w_j \in \partial \mathcal{R}, \tag{2.10}$$

where $\mathcal{R} = \text{conv}(f(I) \cup f(-I))$ denotes the Elfving set [see Elfving (1952)]. Consequently, if all weights in (2.9) are nonnegative, it follows from Elfving's theorem that the design

$$\xi_{c^*}^* = \begin{pmatrix} s_1 & \cdots & s_m \\ w_1 & \cdots & w_m \end{pmatrix} \tag{2.11}$$

is $c^*$-optimal in the regression model (2.8) [see Elfving (1952)], where $c^* \in \mathbb{R}^m$ denotes the vector of coefficients of the Chebyshev polynomial defined in the previous paragraph. The following results relate this design to the $E$-optimal design.

LEMMA 2.1. *Assume that $f_1, \ldots, f_m$ generate a Chebyshev system on the interval $I$ such that the Chebyshev points are unique. If the minimum eigenvalue of the information matrix of an E-optimal design has multiplicity one, then the design $\xi_{c^*}^*$ defined by (2.9) and (2.11) is E-optimal in the regression model (2.8). Moreover, in this case the E-optimal design is unique.*

LEMMA 2.2. *Assume that the functions $f_1, \ldots, f_m$ generate a Chebyshev system on the interval $I$ with Chebyshev polynomial $c^{*T}f(t)$ and let $\xi_{c^*}^*$ denote the $c^*$-optimal design in the regression model (2.2) defined by (2.11). Then $c^*$ is an eigenvector of the information matrix $M(\xi_{c^*}^*, b)$, and if the corresponding eigenvalue $\lambda = \frac{1}{\|c^*\|^2}$ is the minimal eigenvalue, then $\xi_{c^*}^*$ is also E-optimal in the regression model (2.8).*

We now discuss the $c$-optimal design problem in the regression model (2.8) for a general vector $c \in \mathbb{R}^m$ (not necessarily equal to the vector $c^*$ of coefficients of the Chebyshev polynomial). Assume again that $f_1, \ldots, f_m$ generate a Chebyshev system on the interval $I$. As a candidate for the $c$-optimal design, we consider the measure

$$(2.12) \qquad \xi_c = \xi_c(b) = \begin{pmatrix} s_1 & \cdots & s_m \\ w_1 & \cdots & w_m \end{pmatrix},$$

where the support points are the Chebyshev points and the weights are already chosen such that the expression $c^T M^{-1}(\xi_c, b) c$ becomes minimal, that is,

$$(2.13) \qquad w_i = \frac{|e_i^T J F^{-1} c|}{\sum_{j=1}^m |e_j^T J F^{-1} c|}, \qquad i = 1, \ldots, m,$$

where $e_j = (0, \ldots, 0, 1, 0, \ldots, 0)^T \in \mathbb{R}^m$ denotes the $j$th unit vector [see Kitsos, Titterington and Torsney (1988), Pukelsheim and Torsney (1991) or Pukelsheim (1993)]. The following result characterizes the optimal designs for estimating the individual coefficients.

LEMMA 2.3. *Assume that the functions $f_1, \ldots, f_m$ generate a Chebyshev system on the interval $I$. The design $\xi_{e_j}$ defined by (2.12) and (2.13) for the vector $c = e_j$ is $e_j$-optimal in the linear regression model (2.8) if the system $\{f_i | i \in \{1, \ldots, m\} \setminus \{j\}\}$ is a weak Chebyshev system on the interval $I$.*

If the sufficient conditions of Lemma 2.3 are not satisfied, the determination of the $e_j$-optimal designs is a substantially harder problem and optimal designs for estimating individual coefficients have only been found numerically in rare circumstances [see Studden (1968) or Dette, Melas and Pepelyshev (2004)]. In many cases the resulting designs yield a singular information matrix, which makes its determination by standard methods difficult.

REMARK 2.4. It is worthwhile to mention that, in general, the sufficient condition of Lemma 2.3 is not satisfied in the regression model (2.2). To see this, assume that $k \geq 3$, that the function $\varphi$ is continuously differentiable with respect to the second argument and that the functions $f_1(\cdot, b), \ldots, f_m(\cdot, b)$ defined by (2.1) generate a Chebyshev system for any $b$. Define an $(m-1) \times (m-1)$ matrix

$$F_j(x) := (h_1(t_i), \ldots, h_s(t_i), \varphi(t_i, b_1), \varphi'(t_i, b_1), \ldots, \varphi(t_i, b_{j-1}),$$
$$\varphi'(t_i, b_{j-1}), \varphi(t_i, x), \varphi(t_i, b_{j+1}), \ldots, \varphi(t_i, b_k), \varphi'(t_i, b_k))_{i=1}^{m-1},$$

where $c < t_1 < \cdots < t_{m-1} < d$, $b_i \neq b_j$ whenever $i \neq j$ and $x \neq b_i$. We choose $t_1, \ldots, t_{m-1}$ such that $g(x) = \det F_j(x) \not\equiv 0$ (note that the functions $f_1, \ldots, f_m$ form a Chebyshev system and, therefore, this is always possible) and observe that $g(b_i) = 0$, $i = 1, \ldots, k; i \neq j$. Because $k \geq 3$ and $g$ is continuously differentiable, it follows that there exist two points, say $x^*$ and $x^{**}$, such that $g'(x^*) < 0$ and $g'(x^{**}) > 0$. Consequently, there exists an $\bar{x}$ such that

$$0 = g'(\bar{x}) = \det (f_\nu(t_i, b_{\bar{x}}))_{i=1,\ldots,m-1}^{\nu=1,\ldots,m, \nu \neq s+2j-1},$$

where the vector $b_{\bar{x}}$ is defined by $b_{\bar{x}} = (b_1, \ldots, b_{j-1}, \bar{x}, b_{j+1}, \ldots, b_k)^T$. Note that the Chebyshev property of the functions $f_1, \ldots, f_{s+2j-2}, f_{s+2j}, \ldots, f_m$ would imply that all determinants in (2.5) are of the same sign (otherwise there exists a $b$ such that the determinant vanishes for $t_1 < \cdots < t_{m-1}$). Therefore, the conditions $g'(x^*) < 0$, $g'(x^{**}) > 0$ yield that there exists an $\tilde{x} \in (x^*, \bar{x})$ or $\tilde{x} \in (\bar{x}, x^{**})$, such that the system of regression functions

$$\{f_1(t, b_{\tilde{x}}), \ldots, f_{s+2j-2}(t, b_{\tilde{x}}), f_{s+2j}(t, b_{\tilde{x}}), \ldots, f_m(t, b_{\tilde{x}})\}$$
$$= \{h_1(t), \ldots, h_s(t), \varphi(t, b_1), \varphi'(t, b_1), \ldots, \varphi'(t, b_{j-1}),$$
$$\varphi'(t, \tilde{x}), \varphi(t, b_{j+1}), \varphi'(t, b_{j+1}), \ldots, \varphi'(t, b_k)\}$$

is not a weak Chebyshev system on the interval $I$. Finally, in the case $k = 2$, if $\lim_{|b| \to \infty} \varphi(t, b) \to 0$, it can be shown by a similar argument that there exists an $\tilde{x}$ such that the system

$$\{h_1(t), \ldots, h_s(t), \varphi(t, b_1), \varphi'(t, b_1), \varphi'(t, \tilde{x})\}$$

is not a Chebyshev system on the interval $I$.

**3. Asymptotic analysis of $E$- and $c$-optimal designs.** Recall the definition of the information matrix in (2.3) for the model (2.2) with design space given by $I = [c_1, d_1]$ and assume that the nonlinear parameters vary in a compact interval, say $b_i \in [c_2, d_2]$, $i = 1, \ldots, k$. We are interested in the asymptotic properties of $E$- and $c$-optimal designs if

(3.1) $$b_i = x + \delta r_i, \qquad i = 1, \ldots, k,$$

for some fixed $x \in [c_2, d_2]$, fixed $r_1 < r_2 < \cdots < r_k$ and positive $\delta$ satisfying $\delta \to 0$. Note that condition (3.1) implies that all parameters $b_i$ converge to $x$ at the same rate $\delta$. For the asymptotic investigations we study for fixed $\varepsilon, \Delta > 0$ the set

$$
\begin{aligned}
\Omega_{\varepsilon,\Delta} = \Big\{ b \in \mathbb{R}^k | b_i - b_j = \delta(r_i - r_j); \\
i, j = 1, \ldots, k;\ \delta \leq \varepsilon;\ b_i \in [c_2, d_2],\ \min_{i \neq j} |r_i - r_j| \geq \Delta \Big\},
\end{aligned}
\tag{3.2}
$$

introduce the functions

$$
\begin{aligned}
\bar{f}_i(t, x) &= \bar{f}_i(t) = h_i(t), & i &= 1, \ldots, s, \\
\bar{f}_{s+i}(t, x) &= \bar{f}_{s+i}(t) = \varphi^{(i-1)}(t, x), & i &= 1, \ldots, 2k,
\end{aligned}
\tag{3.3}
$$

and the corresponding vector of regression functions

$$
\bar{f}(t, x) = (\bar{f}_1(t, x), \ldots, \bar{f}_{s+2k}(t, x))^T,
\tag{3.4}
$$

where the derivatives of $\varphi(t, x)$ are taken with respect to the second argument. Again the dependency of the functions $\bar{f}_i$ on the parameter $x$ will be omitted whenever it is clear from the context. Note that for a sufficiently smooth function $\varphi$, a simple Taylor expansion shows that under assumption (3.1),

$$
\begin{aligned}
&(\varphi(t, b_1), \varphi'(t, b_1), \ldots, \varphi(t, b_k))^T \\
&\quad = Q(\varphi(t, x), \varphi'(t, x), \ldots, \varphi^{(2k-1)}(t, x))^T + o(1) \\
&\quad = Q\bar{f}(t, x) + o(1)
\end{aligned}
$$

for an appropriately defined matrix $Q \in \mathbb{R}^{2k \times 2k}$ (see the proof of Theorem B.1 in the Appendix). Therefore, optimal designs in the linear model with vector of regression functions given by (3.4) will serve as an approximation for the optimal design in model (2.2) if the parameters $b_i$ are sufficiently close in the sense of (3.1). The following results make this statement more precise.

LEMMA 3.1. *Assume that the function $\varphi : [c_1, d_1] \times [c_2, d_2] \to \mathbb{R}$ in model (1.1) satisfies*

$$
\varphi \in C^{0, 2k-1}([c_1, d_1] \times [c_2, d_2])
$$

*and that for any fixed $x \in [c_2, d_2]$, the functions $\bar{f}_1, \ldots, \bar{f}_{s+2k}$ defined by (3.3) form a Chebyshev system on the interval $[c_1, d_1]$. For any $\Delta > 0$ and any design on the interval $[c_1, d_1]$ with at least $m = s + 2k$ support points, there exists an $\varepsilon > 0$ such that for all $b \in \Omega_{\varepsilon,\Delta}$, the maximum eigenvalue of the inverse information matrix $M^{-1}(\xi, b)$ defined in (2.3) is simple.*

THEOREM 3.2. *Assume that the function $\varphi : [c_1, d_1] \times [c_2, d_2] \to \mathbb{R}$ in model* (1.1) *satisfies*

$$\varphi \in C^{0,2k-1}([c_1, d_1] \times [c_2, d_2])$$

*and that the systems of functions $\{f_1(t,b), \ldots, f_m(t,b)\}$ and $\{\bar{f}_1(t,x), \ldots, \bar{f}_m(t,x)\}$ defined by* (2.1) *and* (3.3), *respectively, are Chebyshev systems on the interval $[c_1, d_1]$ (for arbitrary but fixed $b_1, \ldots, b_k$, $x \in [c_2, d_2]$ with $b_i \neq b_j$ whenever $i \neq j$). If $\varepsilon$ is sufficiently small, then for any $b \in \Omega_{\varepsilon, \Delta}$, the design $\xi^*_{c^*}$ defined by* (2.9) *and* (2.11) *is the unique $E$-optimal design in the regression model* (2.2).

Note that for $b \in \Omega_{\varepsilon, \Delta}$, the $E$-optimal designs can be obtained explicitly by Lemmas 2.1 and 2.2. The support points are the extremal points of the Chebyshev polynomial corresponding to the functions in (2.1), while the weights are given by (2.9).

From Remark 2.4 we may expect that, in general, $c$-optimal designs in the regression model (1.1) are not necessarily supported at the Chebyshev points. Nevertheless, an analogue of Lemma 3.1 is available for specific vectors $c \in \mathbb{R}^m$. The proof is similar to the proof of Lemma 3.1 and, therefore, omitted.

LEMMA 3.3. *Let $e_i = (0, \ldots, 0, 1, 0, \ldots, 0)^T$ denote the $i$th unit vector in $\mathbb{R}^m$. Under the assumptions of Lemma* 3.1 *define a vector $\tilde{\gamma} = (0, \ldots, 0, \gamma_1, \ldots, \gamma_{2k}) \in \mathbb{R}^m$ by*

$$(3.5) \quad \gamma_{2i} = \prod_{j \neq i}(r_i - r_j)^{-2}, \qquad \gamma_{2i-1} = -\gamma_{2i} \sum_{j \neq i} \frac{2}{r_i - r_j}, \qquad i = 1, \ldots, k.$$

(i) *If $c \in \mathbb{R}^m$ satisfies $c^T \tilde{\gamma} \neq 0$, then for any $\Delta > 0$, sufficiently small $\varepsilon$ and any $b \in \Omega_{\varepsilon, \Delta}$, the design $\xi_c(b)$ defined in* (2.12) *and* (2.13) *is $c$-optimal in the regression model* (2.2).

(ii) *The assumption $c^T \tilde{\gamma} \neq 0$ is, in particular, satisfied for the vector $c = e_{s+2j-1}$ for any $j = 1, \ldots, k$ and for the vector $c = e_{s+2j}$ for any $j = 1, \ldots, k$, which satisfies the condition*

$$(3.6) \quad \sum_{\ell \neq j} \frac{1}{r_j - r_\ell} \neq 0.$$

REMARK 3.4. As pointed out by a referee, some explanation of the set $\Omega_{\varepsilon, \Delta}$ is helpful at this point.

Note that the quantity $\Delta \leq \min_{i \neq j} |r_i - r_j|$ yields some mild restriction for the $r_i$ in (3.1) and $\varepsilon$ can be considered as a cut-off point, such that whenever $\delta < \varepsilon$ in (3.1), the statements of Theorem 3.2 and Lemmas 3.1 and 3.3 apply to the corresponding vector $b \in \Omega_{\varepsilon,\Delta}$. This cut-off point cannot be determined explicitly because it depends in a complicated way on $\Delta$, the intervals $[c_1, d_1], [c_2, d_2]$ and the basic function $\varphi(t,x)$ used in the regression model (1.1). Roughly speaking, the results of Lemmas 3.1 and 3.3 and Theorem 3.2 hold for any vector $b$ in a compact neighborhood of the vector $(x, \ldots, x) \in \mathbb{R}^{2k}$. In the examples for the rational model discussed in Section 4 the set $\Omega_{\varepsilon,\Delta}$ coincides with the set of all admissible values for parameter $b$.

Note also that it follows from the proof of Lemma 3.1 that the assumption of compactness of the intervals $[c_1, d_1]$ and $[c_2, d_2]$ is only required for the existence of the set $\Omega_{\varepsilon,\Delta}$. In other words, if condition (3.1) is satisfied and $\delta$ is sufficiently small, the maximum eigenvalue of the matrix $M^{-1}(\xi, b)$ will have multiplicity one (independently of the domain of the function $\varphi$). The same remark applies to the statements of Theorem 3.2 and Lemma 3.3.

Our final result of this section shows that under assumption (3.1) with small $\delta$, the local $E$- and local $c$-optimal designs for the vectors $c$ considered in Lemma 3.3 of Remark 3.4 are very close. To be precise, we assume that the assumptions of Theorem 3.2 are valid and consider the design

$$(3.7) \qquad \bar{\xi}_c = \bar{\xi}_c(x) = \begin{pmatrix} \bar{s}_1 & \cdots & \bar{s}_m \\ \bar{w}_1 & \cdots & \bar{w}_m \end{pmatrix},$$

where $\bar{s}_1, \ldots, \bar{s}_m$ are the Chebyshev points corresponding to the system $\{\bar{f}_i | i = 1, \ldots, m\}$ defined in (3.3),

$$(3.8) \qquad \bar{w}_i = \frac{|e_i^T J \overline{F}^{-1} c|}{\sum_{j=1}^m |e_j^T J \overline{F}^{-1} c|}, \qquad i = 1, \ldots, m,$$

with $\overline{F} = (f_i(\bar{s}_j))_{i,j=1}^m$ and $c \in \mathbb{R}^m$ a fixed vector.

THEOREM 3.5. *Assume that the assumptions of Theorem 3.2 are satisfied and that for the system $\{\bar{f}_1, \ldots, \bar{f}_m\}$ the Chebyshev points are unique.*

(i) *If $\delta \to 0$, the design $\xi_{c^*}^*(b)$ defined by (2.11) and (2.9) converges weakly to the design $\bar{\xi}_{e_m}(x)$ defined by (3.7) and (3.8) for $c = e_m$.*

(ii) *If $c \in \mathbb{R}^m$ satisfies $c^T \tilde{\gamma} \neq 0$ for the vector $\tilde{\gamma}$ with components defined in (3.5) and $\delta \to 0$, then the design $\xi_c^*(b)$ defined by (2.12) and (2.13) converges weakly to the design $\bar{\xi}_{e_m}(x)$.*

(iii) *The assumption $c^T \tilde{\gamma} \neq 0$ is, in particular, satisfied for the vector $c = e_{s+2j-1}$ for any $j = 1, \ldots, k$ and for the vector $c = e_{s+2j}$ for any $j = 1, \ldots, k$, which satisfies condition (3.6).*

REMARK 3.6. Note that Theorem 3.2, Lemma 3.3 and Theorem 3.5 remain valid for the local optimal designs in the nonlinear regression model (1.1).

This follows by a careful inspection of the proofs of the previous results. For example, there exists a set $\Omega_{\varepsilon,\Delta}$ such that for all $b \in \Omega_{\varepsilon,\Delta}$, the maximum eigenvalue of the inverse information matrix in the model (1.1) is simple. Similarly, if $\delta \to 0$ and (3.1) is satisfied, $c$-optimal designs in the nonlinear regression model are given by the design $\xi_{\bar{c}}(b)$ in (2.12) and (2.13) with $\bar{c} = K_a c$, whenever $\tilde{\gamma}^T \bar{c} \neq 0$, and all these designs converge weakly to the $e_m$-optimal design in the linear regression model defined by the functions (3.4).

We finally remark that Theorem 3.5 and Remark 3.6 indicate that $E$-optimal designs are very efficient for estimating the parameters $a_{s+1}, b_1, \ldots, a_{s+k}, b_k$ in the nonlinear regression model (1.1) and the linear model (2.2), because for small differences $|b_i - b_j|$ the $E$-optimal design and the optimal design for estimating an individual coefficient $b_i$ $(i = 1, \ldots, k)$ are close to the optimal design for estimating the coefficient $b_k$. We will illustrate this fact in the following section, which discusses the rational model in more detail.

**4. Rational models.** In this section we discuss the rational model (1.2) in more detail, where the design space is a compact or seminfinite interval $I$. In contrast to the work of Imhof and Studden (2001), we assume that the nonlinear parameters $b_1, \ldots, b_k \notin I$ are not known by the experimenter but have to be estimated from the data. A typical application of this model can be found in the work of Dudzinski and Mykytowycz (1961), where this model was used to describe the relation between the weight of the dried eye lens of the European rabbit and the age of the animal. In the notation of Sections 2 and 3 we have $f(t) = f(t, b) = (f_1(t), \ldots, f_m(t))^T$, with

$$f_i(t) = t^{i-1}, \qquad i = 1, \ldots, s,$$

(4.1)  $$f_{s+2i-1}(t) = f_{s+2i-1}(t, b) = \frac{1}{t - b_i},$$

$$f_{s+2i}(t) = f_{s+2i}(t, b) = \frac{1}{(t - b_i)^2}, \qquad i = 1, \ldots, k,$$

and the equivalent linear regression model is given by (1.4). The corresponding limiting model is determined by the regression functions $\bar{f}(t) = \bar{f}(t) = (\bar{f}_1(t), \ldots, \bar{f}_m(t))^T$, with

(4.2)  $$\bar{f}_i(t) = t^{i-1}, \qquad \bar{f}_{i+s}(t) = \bar{f}_{s+i}(t, x) = \frac{1}{(t-x)^i}, \qquad i = 1, \ldots, s.$$

Some properties of the functions defined by (4.1) and (4.2) are discussed in the following lemma.

LEMMA 4.1. *Define* $\mathcal{B} = \{b = (b_1, \ldots, b_k)^T \in \mathbb{R}^k | b_i \notin I; b_i \neq b_j\}$. *Then the following assertions are true:*

(i) *If $I$ is a finite interval or $I \subset [0, \infty)$ and $b \in \mathcal{B}$, then the system $\{f_1(t_1, b), \ldots, f_m(t, b)\}$ defined in (4.1) is a Chebyshev system on the interval $I$. If $x \notin I$, then the system $\{\bar{f}_1(t, x), \ldots, \bar{f}_m(t, x)\}$ defined by (4.2) is a Chebyshev system on the interval $I$.*

(ii) *Assume that $b \in \mathcal{B}$ and that one of the following conditions is satisfied:*

    (a) *$I \subset [0, \infty)$,*
    (b) *$s = 1$ or $s = 0$.*

*For any $j \in \{1, \ldots, k\}$, the system of regression functions $\{f_i(t, b) | i = 1, \ldots, m, i \neq s + 2j\}$ is a Chebyshev system on the interval $I$.*

(iii) *If $I$ is a finite interval or $I \subset [0, \infty)$, $k \geq 2$, and $j \in \{1, \ldots, k\}$, then there exists a nonempty set $W_j \subset \mathcal{B}$ such that for all $b \in W_j$, the system of functions $\{f_i(t, b) | i = 1, \ldots, m; i \neq s + 2j - 1\}$ is not a Chebyshev system on the interval $I$.*

The case $k = 1$ will be studied more explicitly in Example 4.5. Note that the third part of Lemma 4.1 shows that for $k \geq 2$, the main condition of Theorem 2.1 in the paper of Imhof and Studden (2001) is *not* satisfied in general for the linear regression model with the functions given by (4.1). These authors assumed that every subsystem of $\{f_1, \ldots, f_m\}$ which consists of $m - 1$ of these functions is a weak Chebyshev system on the interval $I$. Because the design problem for this model is equivalent to the design problem for the model (1.2) (where the nonlinear parameters are not known and have to be estimated), it follows that, in general, we cannot expect local $E$-optimal designs for the rational model to be supported at the Chebyshev points. However, the linearized regression model (1.4) is a special case of the general model (2.2) with $\varphi(t, b) = (t - b)^{-1}$ and all results of Section 3 are applicable here. In particular, we obtain that the $E$-optimal designs and the optimal designs for estimating the individual coefficients $a_{s+1}, b_1, \ldots, a_{s+k}, b_k$ are supported at the Chebyshev points if the nonlinear parameters $b_1, \ldots, b_k$ are sufficiently close (see Theorem 3.2, Lemma 3.3 and Remark 3.6).

THEOREM 4.2. (i) *If $s = 1$, then the Chebyshev points $s_1 = s_1(b), \ldots, s_m = s_m(b)$ for the system of regression functions in (4.1) on the interval $[-1, 1]$ are given by the zeros of the polynomial*

$$(4.3) \qquad (1 - t^2) \sum_{i=0}^{4k} d_i U_{-2k+s+i-2}(t),$$

*where $U_j(x)$ denotes the $j$th Chebyshev polynomial of the second kind [see Szegö (1975)], $U_{-1}(x) = 0, U_{-n}(x) = -U_{n-2}(x)$ and the factors $d_0, \ldots, d_{4k}$* *are defined as the coefficients of the polynomial*

$$\sum_{i=0}^{4k} d_i t^i = \prod_{i=1}^{k} (t - \tau_i)^4, \tag{4.4}$$

*where*

$$2b_i = \tau_i + \frac{1}{\tau_i}, \qquad i = 1, \ldots, k.$$

(ii) *Let $\Omega_E \subset \mathcal{B}$ denote the set of all $b$ such that an $E$-optimal design for the model* (1.4) *is given by* (2.11) *and* (2.9). *Then $\Omega_E \neq \varnothing$.*

REMARK 4.3. (a) The Chebyshev points for the system (4.1) on an arbitrary finite interval $I \subset \mathbb{R}$ can be obtained by rescaling the points onto the interval $[-1, 1]$. The case $s = 0$ and $I = [0, \infty)$ will be discussed in more detail in Examples 4.5 and 4.6.

(b) It follows from Theorem 3.2 that the set $\Omega_E$ defined in the second part of Theorem 4.1 contains the set $\Omega_{\varepsilon, \Delta}$ defined in (3.2) for sufficiently small $\varepsilon$. In other words, if the nonlinear parameters $b_1, \ldots, b_k$ are sufficiently close, the local $E$-optimal design will be supported at the Chebyshev points with weights given by (2.9). Moreover, we will demonstrate in the subsequent examples that in many cases the set $\Omega_E$ coincides with the full set $\mathcal{B}$.

(c) In applications the Chebyshev points can be calculated numerically with the Remez algorithm [see Studden and Tsay (1976) or DeVore and Lorentz (1993)]. In some cases these points can be obtained explicitly.

REMARK 4.4. We note that a similar result is valid for $c$-optimal designs in the rational regression model (1.4). For example, assume that one of the assertions of Lemma 4.1 is valid and that we are interested in estimating a linear combination $c^T \beta$ of the parameters in the rational model (1.4). We obtain from Lemma 3.3 that if $c \in \mathbb{R}^m$ satisfies $c^T \tilde{\gamma} \neq 0$, then for sufficiently small $\varepsilon$ and any $b \in \Omega_{\varepsilon, \Delta}$, the design $\xi_c(b)$ defined in (2.12) and (2.13) is $c$-optimal. In particular, this is true for $c = e_{s+2j-1}$ (for all $j = 1, \ldots, k$) and the vector $c = e_{s+2j}$ if the index $j$ satisfies the condition (3.6). Note that due to the third part of Lemma 4.1 in the case $k \geq 2$, there exists a $b \in \mathcal{B}$ such that the $e_{s+2j}$-optimal design is not necessarily supported at the Chebyshev points. However, from Theorem 3.5 it follows that for a vector $b \in \mathcal{B}$ satisfying (3.1) with $\delta \to 0$ and any vector $c$ with $c^T \tilde{\gamma} \neq 0$, we have for the designs $\xi_{c^*}^*(b)$ and $\xi_c^*(b)$ defined by (2.11) and (2.12)

$$\xi_{c^*}^*(b) \to \bar{\xi}_{e_m}(x), \qquad \xi_c^*(b) \to \bar{\xi}_{e_m}(x),$$

where the design $\bar{\xi}_{e_m}(x)$ is defined in (3.7) and (3.8) and is $e_m$-optimal in the limiting model with the regression functions (4.2). We conclude this section with two examples. Further examples considering a finite interval as design space and a comparison with $D$-optimal designs can be found in the technical report of Dette, Melas and Pepelyshev (2002).

EXAMPLE 4.5. Consider the rational model

$$Y = \frac{a}{t-b} + \varepsilon, \qquad t \in [0, \infty), \tag{4.5}$$

with $b < 0$ (here we have $k = 1$, $s = 0$, $I = [0, \infty)$). The corresponding equivalent linear regression model is given by

$$Y = \beta^T f(t, b) + \varepsilon = \frac{\beta_1}{t-b} + \frac{\beta_2}{(t-b)^2} + \varepsilon. \tag{4.6}$$

In this case it follows from the first part of Lemma 4.1 that the system of regression functions $\{\frac{1}{t-b}, \frac{1}{(t-b)^2}\} = \{f_1(t), f_2(t)\}$ is a Chebyshev system on the interval $[0, \infty)$ whenever $b < 0$. Moreover, any subsystem (consisting of one function) is obviously a Chebyshev system on the interval $[0, \infty)$. The Chebyshev points are given by $s_1 = 0$ and $s_2 = \sqrt{2}|b| = -\sqrt{2}b$. Now we consider the design $\xi_c^*(b)$ defined in (2.12) as a candidate for the $c$-optimal design in model (4.6). The weights (for any $c \in \mathbb{R}^2$) are obtained from formula (2.13) and a straightforward calculation shows that the $c$-optimal design $\xi_c^*(b)$ has masses $\omega_1$ and $1 - \omega_1$ at the points $0$ and $\sqrt{2}|b|$, respectively, where

$$\omega_1 = \frac{|b(-\sqrt{2}c_1 + (2+\sqrt{2})c_2 b)|}{|b|\{|-\sqrt{2}c_1 + (2+\sqrt{2})c_2 b| + (4 + 3\sqrt{2})|-c_1 + c_2 b|\}}.$$

It can easily be checked by Elfving's theorem [see Elfving (1952)] or by the equivalence theorem for $c$-optimality [see Pukelsheim (1993)] that this design is, in fact, $c$-optimal in the regression model (4.6) whenever $\frac{c_2}{c_1} \notin [\frac{1}{b}, \frac{1}{(1+\sqrt{2})b}]$. In the remaining cases the $c$-optimal design is a one-point design supported at $t = b - \frac{c_1}{c_2}$. In particular, by Lemma 2.3, the $e_1$- and $e_2$-optimal designs for estimating the coefficients $\beta_1$ and $\beta_2$ in the model (4.6) have weights $\frac{1}{4}(2 - \sqrt{2}), \frac{1}{4}(2 + \sqrt{2})$ and $1 - \frac{1}{\sqrt{2}}, \frac{1}{\sqrt{2}}$ at the points $0, \sqrt{2}|b|$, respectively. It follows from the results of Imhof and Studden (2001) that an $E$-optimal design in the regression model (4.6) is given by the $c^*$-optimal design for the Chebyshev vector $c^* = (1 + \sqrt{2})|b|(-2, |b|(1+\sqrt{2}))^T$, which has masses $w_1$ and $1 - w_1$ at the points $0$ and $\sqrt{2}|b|$, respectively, where

$$w_1 = \frac{1}{2} \frac{(2-\sqrt{2})(6 - 4\sqrt{2} + b^2)}{b^2 + 12 - 8\sqrt{2}} = 1 - \frac{1}{2} \frac{\sqrt{2}(2\sqrt{2} - 2 + b^2)}{b^2 + 12 - 8\sqrt{2}}.$$

Alternatively, the $E$-optimal design could be also obtained by the geometric method of Dette and Haines (1994), which is especially designed for models with two parameters.

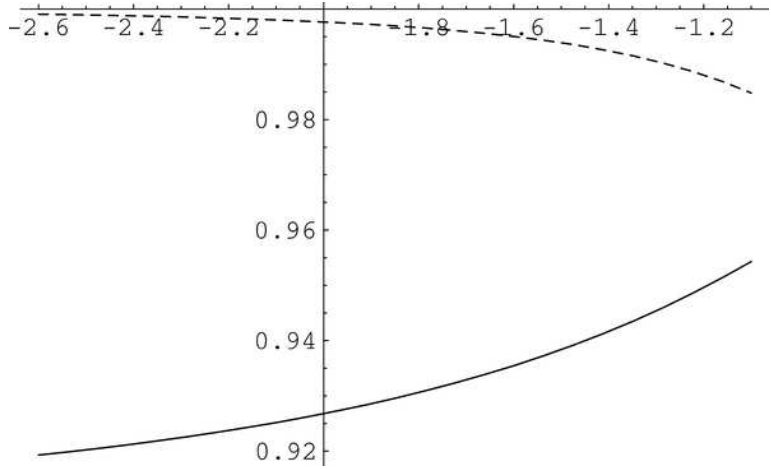

FIG. 1. *Efficiencies of the E-optimal design $\xi^*(b)$ for estimating the individual coefficients in the regression model* (4.6) *for various values of $b \in [-2.5, -1]$. Solid line: $\mathrm{eff}_1(\xi^*(b))$, dotted line: $\mathrm{eff}_2(\xi^*(b))$.*

In Figure 1 we show the efficiencies of the $E$-optimal design for estimating the coefficients $\beta_1$ and $\beta_2$ in the regression model (4.6), that is,

$$\mathrm{eff}_i(\xi_E^*(b))$$
$$= \left( \frac{e_i^T M^{-1}(\xi_E^*(b), b) e_i}{e_i^T M^{-1}(\xi_{e_i}^*, b) e_i} \right)^{-1}$$

(4.7)
$$= \begin{cases} \dfrac{28(b^4(5\sqrt{2} - 7) + b^2(34\sqrt{2} - 48) + 396 - 280\sqrt{2})}{(9\sqrt{2} - 11)(b^2 - 8\sqrt{2} + 12)(7b^2 + 16\sqrt{2} - 20)}, & \text{if } i = 1, \\ \dfrac{b^4(\sqrt{2} - 1) + (6\sqrt{2} - 8)b^2 + 68 - 48\sqrt{2}}{(\sqrt{2} - 1)(b^2 - 8\sqrt{2} + 12)(b^2 - 6\sqrt{2} + 8)}, & \text{if } i = 2 \end{cases}$$

[for technical details for this calculation see Dette, Melas and Pepelyshev (2002)]. We observe for the $e_1$-efficiency for all $b \leq -1$ the inequality

$$0.9061 \approx = \lim_{b \to -\infty} \mathrm{eff}_1(\xi_E^*(b)) \leq \mathrm{eff}_1(\xi_E^*(b)) \leq \mathrm{eff}_1(\xi_E^*(-1)) \approx 0.9595,$$

and similarly for the $e_2$-efficiency

$$0.9805 \approx \mathrm{eff}_2(\xi_E^*(-1)) \leq \mathrm{eff}_2(\xi_E^*(b)) \leq \lim_{b \to -\infty} \mathrm{eff}_2(\xi_E^*(b)) = 1.$$

This demonstrates that the $E$-optimal design yields very accurate estimates for the individual parameters in the regression model (4.6).

We finally mention the results for the local optimal design in the rational model (4.5), which maximize or minimize the corresponding functional for the matrix $K_a^{-1} M(\xi, b) K_a^{-1}$, where $K_a = \mathrm{diag}(1, -\frac{1}{a})$. Obviously, the local $e_1$- and $e_2$-designs coincide with the corresponding designs in the equivalent linear regression model (4.6). On the other hand, the $c$-optimal design for the rational model (4.5) is obtained from the $\bar{c}$-optimal design $\xi_{\bar{c}}^*(b)$ for the model (4.6) with $\bar{c} = K_a c = (c_1, -c_2/a)^T$. Similarly, the local $E$-optimal design for the rational model (4.5) has masses $w_1^*$ and $1 - w_1^*$ at the points $0$ and $\sqrt{2}|b|$, where the weights are given by

$$w_1^* = \frac{2\sqrt{2}a^2 + (4 + 3\sqrt{2})b^2}{2\{4(1+\sqrt{2})a^2 + (7+5\sqrt{2})b^2\}} = 1 - \frac{(4+3\sqrt{2})(2a^2 + (1+\sqrt{2})b^2)}{2\{4(1+\sqrt{2})a^2 + (7+5\sqrt{2})b^2\}}.$$

An investigation of the efficiencies for the $E$-optimal design in the rational model (4.5) yields similar results as in the corresponding equivalent linear regression model (4.6). For a broad range of parameter values $(a, b)$ the local $E$-optimal designs in the rational model (4.5) are very efficient for estimating the individual parameters.

EXAMPLE 4.6. We now discuss $E$-optimal designs for the rational model

$$(4.8) \qquad Y = \frac{a_1}{t - b_1} + \frac{a_2}{t - b_2} + \varepsilon, \qquad t \in [0, \infty),$$

where $b_1, b_2 < 0$; $|b_2 - b_2| > 0$ $(k = 2, s = 0)$. The corresponding equivalent linear regression model is given by

$$(4.9) \qquad Y = \frac{\beta_1}{t - b_1} + \frac{\beta_2}{(t - b_1)^2} + \frac{\beta_3}{t - b_2} + \frac{\beta_4}{(t - b_2)^2} + \varepsilon.$$

The results of Section 3 show that for sufficiently close parameters $b_i$, the $E$- and $e_i$-optimal designs are supported at the Chebyshev points and that the $c^*$-optimal design is the unique $E$-optimal design. In this case the local optimal designs cannot be found explicitly. Therefore, we used these designs for any vector $(b_1, b_2)$ under consideration as candidates for the optimal designs. In other words, we used the Chebyshev points as support points and calculated the optimal weights from the formulas presented in Section 2 to obtain candidates for the local optimal design. The optimality for a concrete choice was finally verified by an application of the results in Section 2 (see the discussion below). For the sake of brevity, we restrict ourselves to model (4.9), which corresponds to the local optimal design problem for model (4.8) with $(a_1, a_2) = (1, 1)$. In our comparison we will also include the $E$-optimal design in the limiting model under assumption (3.1), that is,

$$(4.10) \qquad Y = \frac{\beta_1}{t - x} + \frac{\beta_2}{(t - x)^2} + \frac{\beta_3}{(t - x)^3} + \frac{\beta_4}{(t - x)^4} + \varepsilon,$$

where the parameter $x$ is chosen as $x = (b_1 + b_2)/2$. Without loss of generality we assume that $x = -1$, because in the general case the optimal designs can be obtained by a simple scaling argument. The limiting optimal design was obtained numerically and has masses 0.13, 0.26, 0.27, 0.34 at the points 0, 0.18, 1.08 and 7.9, respectively.

TABLE 1
*E-optimal designs for linear regression model* (4.9) *on the interval* $[0,\infty)$,
*where* $b_1 = -1 - z$, $b_2 = -1 + z$. *These designs are E-optimal in the rational
model* (4.8) *for the initial parameter* $a_1 = a_2 = 1$. *Note that the smallest
support point of the E-optimal design* $(t_{1E}^*)$ *is* 0

| $z$ | 0.1 | 0.2 | 0.3 | 0.4 | 0.5 | 0.6 | 0.7 | 0.8 | 0.9 | 0.95 |
|---|---|---|---|---|---|---|---|---|---|---|
| $t_{2E}^*$ | 0.18 | 0.17 | 0.17 | 0.16 | 0.15 | 0.13 | 0.11 | 0.09 | 0.05 | 0.03 |
| $t_{3E}^*$ | 1.08 | 1.06 | 1.03 | 0.99 | 0.94 | 0.87 | 0.77 | 0.65 | 0.47 | 0.34 |
| $t_{4E}^*$ | 7.85 | 7.77 | 7.65 | 7.46 | 7.21 | 6.88 | 6.45 | 5.88 | 5.05 | 4.43 |
| $w_{1E}^*$ | 0.13 | 0.13 | 0.13 | 0.13 | 0.12 | 0.10 | 0.08 | 0.07 | 0.05 | 0.03 |
| $w_{2E}^*$ | 0.26 | 0.26 | 0.27 | 0.26 | 0.25 | 0.22 | 0.20 | 0.17 | 0.13 | 0.10 |
| $w_{3E}^*$ | 0.27 | 0.27 | 0.28 | 0.28 | 0.28 | 0.28 | 0.28 | 0.28 | 0.28 | 0.28 |
| $w_{4E}^*$ | 0.34 | 0.33 | 0.33 | 0.33 | 0.36 | 0.39 | 0.44 | 0.49 | 0.54 | 0.59 |

Theorem 3.2 shows that for sufficiently small $|b_1 - b_2|$, E-optimal designs for the model (4.9) are given by the design $\xi_{c^*}^*(b)$ defined in (2.9) and (2.11). From Lemma 2.2 it follows that the design $\xi_{c^*}^*(b)$ is E-optimal whenever

$$\lambda_{c^*} := \frac{c^{*T} M(\xi_E^*(b), b) c^*}{c^{*T} c^*} \leq \lambda_{(2)}(M(\xi_E^*(b), b)) = \lambda_{(2)},$$

where $\lambda_{\min}(M(\xi_E^*(b), b)) \leq \lambda_{(2)} \leq \cdots \leq \lambda_{(m)}$ denote the ordered eigenvalues of the matrix $M(\xi_E^*(b), b)$. The ratio $\lambda_{(2)}/\lambda_{c^*}$ is illustratively depicted in Fig-

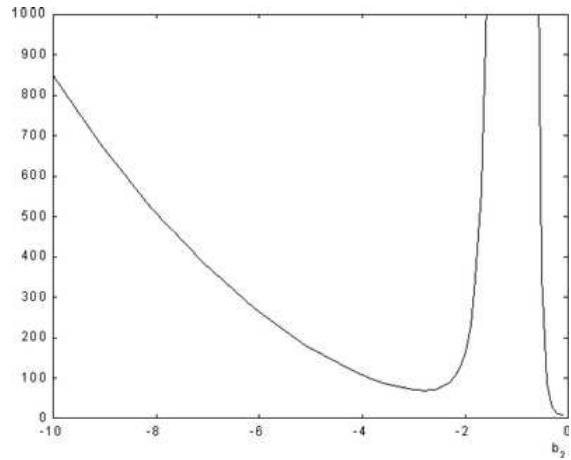

FIG. 2.   *The ratio* $\lambda_{(2)}/\lambda_{c^*}$ *for the design* $\xi_E^*(b)$, *where* $b = (-1, b_2)$. *The designs are E-optimal if this ratio is greater than or equal to* 1.

OPTIMAL DESIGNS 19

ure 2 for $b_1 = 1$ and a broad range of $b_2$ values, which shows that it is always bigger than 1. Other cases yield a similar picture and, in practice, the local $E$-optimal design for the rational model (4.8) and the equivalent linear regression model (4.9) is always supported at the Chebyshev points and given by (2.9) and (2.11). In Tables 1 and 2 we give the main characteristics and efficiencies for the local $E$-optimal design $\xi_E^*(b)$ and for the $E$-optimal design $\bar{\xi}_E^*(\frac{b_1+b_2}{2})$ in the limiting regression model (4.10). The efficiencies are defined by (4.7) and we observe again very good performance of the $E$-optimal designs. The behavior of the design $\bar{\xi}_E$ in the limiting regression model (4.9) is interesting from a practical point of view because it is very similar to the performance of the $E$-optimal design for a broad range of $b_1$ and $b_2$ values. Consequently, this design might be appropriate if rather imprecise prior information for the nonlinear parameters is available. For example, if it is known (from scientific considerations) that $b_1 \in [\underline{b}_1, \bar{b}_1]$, $b_2 \in [\underline{b}_2, \bar{b}_2]$, the design $\bar{\xi}_E(\frac{\underline{b}_1+\bar{b}_2}{2})$ might be a robust choice for practical experiments.

## APPENDIX A: PROOFS

PROOF OF LEMMA 2.1. Let $\xi_E^*$ denote an $E$-optimal design such that the minimum eigenvalue $\lambda = \lambda_{\min}(M(\xi_E^*, b))$ of the information matrix $M(\xi_E^*, b)$ has multiplicity one with corresponding eigenvector $z \in \mathbb{R}^m$. By the equivalence theorem for the $E$-optimality criterion [see Pukelsheim (1993), pages 181 and 182], we obtain for the matrix $E = zz^T/\lambda$,

$$\left(\frac{1}{\sqrt{\lambda}}z^T f(t)\right)^2 = f^T(t)Ef(t) \le 1$$

for all $t \in I$ with equality at the support points of $\xi_E^*$. Because the Chebyshev polynomial is unique it follows that (up to the factor $\mp 1$) $c^* = \frac{1}{\sqrt{\lambda}}z$ and that $\mathrm{supp}(\xi_E^*) = \{s_1, \ldots, s_m\}$. Now Theorem 3.2 in Dette and Studden (1993) implies that $\xi_E^*$ is also $c^*$-optimal, where $c^* \in \mathbb{R}^m$ denotes the vector of coefficients of the Chebyshev polynomial. Consequently, by the discussion of the previous paragraph we have $\xi_E^* = \xi_{c^*}^*$, which proves the assertion. $\square$

TABLE 2
The efficiency (4.7) of the $E$-optimal designs $\xi_E^*$ in the linear regression model (4.9) on the interval $[0, \infty)$ with $b_1 = -1 - z$, $b_2 = -1 + z$ and the efficiency of the $E$-optimal design $\bar{\xi}_E^*(-1)$ in the corresponding limiting model (4.10)

| $z$ | 0.1 | 0.2 | 0.3 | 0.4 | 0.5 | 0.6 | 0.7 | 0.8 | 0.9 | 0.95 |
|---|---|---|---|---|---|---|---|---|---|---|
| $\text{eff}_1(\xi_E^*)$ | 1.00 | 1.00 | 0.99 | 0.94 | 0.70 | 0.45 | 0.50 | 0.55 | 0.64 | 0.78 |
| $\text{eff}_2(\xi_E^*)$ | 0.99 | 0.99 | 0.98 | 0.98 | 0.99 | 1.00 | 1.00 | 1.00 | 1.00 | 1.00 |
| $\text{eff}_3(\xi_E^*)$ | 1.00 | 1.00 | 1.00 | 0.99 | 0.95 | 0.87 | 0.76 | 0.68 | 0.58 | 0.44 |
| $\text{eff}_4(\xi_E^*)$ | 1.00 | 0.99 | 0.98 | 0.94 | 0.87 | 0.76 | 0.62 | 0.54 | 0.44 | 0.31 |
| $\text{eff}_1(\bar{\xi}_E^*(-1))$ | 1.00 | 0.99 | 0.98 | 0.94 | 0.79 | 0.61 | 0.39 | 0.32 | 0.29 | 0.27 |
| $\text{eff}_2(\bar{\xi}_E^*(-1))$ | 0.99 | 0.97 | 0.94 | 0.88 | 0.78 | 0.65 | 0.49 | 0.40 | 0.31 | 0.21 |
| $\text{eff}_3(\bar{\xi}_E^*(-1))$ | 1.00 | 0.99 | 0.98 | 0.95 | 0.88 | 0.75 | 0.54 | 0.40 | 0.24 | 0.08 |
| $\text{eff}_4(\bar{\xi}_E^*(-1))$ | 1.00 | 0.99 | 0.98 | 0.95 | 0.90 | 0.78 | 0.57 | 0.41 | 0.24 | 0.07 |

PROOF OF LEMMA 2.2. From the identity (2.10) and the Chebyshev property (2.7) it follows immediately that $c^*$ is an eigenvector of the matrix

$$M(\xi_{c^*}^*, b) = \sum_{i=1}^m f(s_i) f^T(s_i) w_i$$

with corresponding eigenvalue $\lambda = 1/\|c^*\|^2$. Now if $\lambda = \lambda_{\min}(M(\xi_{c^*}^*, b))$, we define the matrix $E = \lambda c^* c^{*T}$ and obtain from the Chebyshev properties (2.7) that

$$f^T(t) E f(t) = \lambda (c^{*T} f(t))^2 \leq \lambda = \lambda_{\min}(M(\xi_{c^*}^*, b))$$

for all $t \in I$. The assertion of Lemma 2.2 now follows from the equivalence theorem for $E$-optimality [see Pukelsheim (1993)]. $\square$

PROOF OF LEMMA 2.3. If $f_1, \ldots, f_m$ generate a weak Chebyshev system on the interval $I$, it follows from Theorem 2.1 in Studden (1968) that the design $\xi_{e_j}$ defined in (2.12) and (2.13) is $e_j$-optimal if

$$\varepsilon e_i^T J F^{-1} e_j \geq 0, \qquad i = 1, \ldots, m,$$

for some $\varepsilon \in \{-1, 1\}$. The assertion of Lemma 2.3 is now obtained by Cramér's rule. $\square$

PROOF OF LEMMA 3.1. Recall the definition of the functions in (3.3) and let

$$(\text{A.1}) \qquad \overline{M}(\xi, x) = \int_c^d \bar{f}(t, x) \bar{f}^T(t, x) \, d\xi(x)$$

denote the information matrix in the corresponding linear regression model. Because of the Chebyshev property of the functions $\bar{f}_1, \ldots, \bar{f}_{s+2k}$, we have $|\overline{M}(\xi, x)| \neq 0$ (note that the design $\xi$ has at least $s + 2k$ support points). It follows from Theorem B.1 that under the condition (3.1) with $\delta \to 0$, the asymptotic expansion

$$(\text{A.2}) \qquad \delta^{4k-2} M^{-1}(\xi, b) = h \bar{\gamma} \bar{\gamma}^T + o(1)$$

is valid, where the vector $\bar{\gamma} = (\bar{\gamma}_1, \ldots, \bar{\gamma}_{s+2k})^T$ is defined by

$$\bar{\gamma}_{s+2i-1} = -\prod_{j \neq i}(r_i - r_j)^{-2} \cdot \sum_{j \neq i} \frac{2}{r_i - r_j}, \qquad i = 1, \ldots, k,$$

(A.3)
$$\bar{\gamma}_1 = \cdots = \bar{\gamma}_s = 0, \qquad \bar{\gamma}_{s+2i} = 0, \qquad i = 1, \ldots, k,$$

and the constant $h$ is given by

(A.4) $$h = ((2k-1)!)^2 (\overline{M}^{-1}(\xi, x))_{m,m}.$$

From (A.2) we obtain that the maximal eigenvalue of the matrix $M^{-1}(\xi, b)$ is simple if $\delta$ is sufficiently small.

For a fixed value $r = (r_1, \ldots, r_k)$ and fixed $x \in \mathbb{R}$ in the representation (3.1), denote by $\varepsilon = \varepsilon(x, r)$ the maximal value (possibly $\infty$) such that the matrix $M^{-1}(\xi, b)$ has a simple maximal eigenvalue for all $\delta \leq \varepsilon$. Then the function $\varepsilon : (x, r) \to \varepsilon(x, r)$ is continuous and the infimum

$$\inf \left\{ \varepsilon(x, b) \Big| x \in [c_1, d_1], \min_{i \neq j} |r_i - r_j| \geq \Delta, \|r\|_2 = 1 \right\}$$

is attained for some $x^* \in [c_1, d_1]$ and $r^*$, which implies $\varepsilon^* = \varepsilon(x^*, r^*) > 0$. This means that for any $b \in \Omega_{\varepsilon^*, \Delta}$, the multiplicity of the maximal eigenvalue of the information matrix $M^{-1}(\xi, b)$ is equal to one. $\square$

PROOF OF THEOREM 3.2. The proof is a direct consequence of Lemma 2.2 and Lemma 3.1, which shows that the multiplicity of the maximum eigenvalue of the inverse information matrix of any design has multiplicity one, if $b \in \Omega_{\varepsilon, \Delta}$ and $\varepsilon$ is sufficiently small. $\square$

PROOF OF THEOREM 3.5. It follows from Theorem 3.2 that the design $\xi_{c^*}^* = \xi_{c^*}^*(b)$ is local $E$-optimal for sufficiently small $\delta > 0$. In other words, if $\delta$ is sufficiently small, the design $\xi_{c^*}^*$ minimizes $\max_{\|c\|_2 = 1} c^T M^{-1}(\xi, b) c$ in the class of all designs. Note that the components of the vector $r = (r_1, \ldots, r_k)$ are ordered, which implies

$$e_{s+2i-1}^T \tilde{\gamma} \neq 0, \qquad i = 1, \ldots, k.$$

Multiplying (B.1) by $\delta^{4k-2}$, it then follows from Theorem B.1 that for some subsequence $\delta_k \to 0 : \xi_{c^*}^* \to \hat{\xi}(x)$, where the design $\hat{\xi}(x)$ minimizes the function

$$\max_{\|c\|_2 = 1} (c^T \tilde{\gamma})^2 e_m^T \overline{M}^{-1}(\xi, x) e_m$$

and the vector $\tilde{\gamma}$ is defined by (A.3). The maximum is attained for $c = \tilde{\gamma}/\|\tilde{\gamma}\|_2$ (independently of the design $\xi$) and, consequently, $\hat{\xi}(x)$ is $e_m$-optimal in the linear regression model defined by the vector of regression functions in (3.4). Now the functions $\bar{f}_1, \ldots, \bar{f}_m$ generate a Chebyshev system and the corresponding Chebyshev points are unique, which implies that the $e_m$-optimal design $\bar{\xi}_{e_m}(x)$ is unique. Consequently, every subsequence of designs $\xi_{c^*}^*(b)$ contains a weakly convergent subsequence with limit $\bar{\xi}_{e_m}(x)$ and this proves the first part of the assertion. For a proof of the second part we note

that a $c$-optimal design minimizes $c^T M^{-1}(\xi,b)c$ in the class of all designs on the interval $I$. Now if $c^T \tilde{\gamma} \neq 0$ and

$$e_{s+2i-1}^T \tilde{\gamma} = -\prod_{j \neq i}(r_i - r_j)^{-2} \sum_{j \neq i} \frac{2}{r_i - r_j} \neq 0$$

for some $i = 1, \ldots, k$, the same argument as in the previous paragraph shows that $\xi_c^*(b)$ converges weakly to the design which maximizes the function $(\tilde{\gamma}^T c)^2 e_m^T \overline{M}^{-1}(\xi, x) e_m$. If $e_{s+2i-1}^T \tilde{\gamma} = 0$ for all $i = 1, \ldots, k$, the condition $c^T \tilde{\gamma} \neq 0$ implies $e_{s+2i}^T \tilde{\gamma} \neq 0$ for some $i = 1, \ldots, k$ and the assertion follows by multiplying (B.1) by $\delta^{4k-4}$ and similar arguments. Finally, the third assertion follows directly from the definition of the vector $\tilde{\gamma}$ in (3.5).  □

PROOF OF LEMMA 4.1.　Part (iii) follows from Remark 2.4. Parts (i) and (ii) are proved similarly and we restrict ourselves to the first case. For this purpose we introduce the functions $\psi(t,b) = (\psi_1(t,\tilde{b}), \ldots, \psi_m(t,\tilde{b}))^T$ with

(A.5)
$$\psi_i(t,\tilde{b}) = t^{i-1}, \qquad i = 1, \ldots, s,$$
$$\psi_{s+i}(t,\tilde{b}) = \frac{1}{t - \tilde{b}_i}, \qquad i = 1, \ldots, 2k,$$

where $\tilde{b} = (\tilde{b}_1, \ldots, \tilde{b}_{2k})^T$ is a fixed vector with $\tilde{b}_i \neq \tilde{b}_j$ if $i \neq j$. With the notation

$$L(\Delta) = \begin{pmatrix} I_s & 0 \\ 0 & G_k(\Delta) \end{pmatrix} \in \mathbb{R}^{m \times m},$$

$$G_k(\Delta) = \begin{pmatrix} G(\Delta) & & \\ & \ddots & \\ & & G(\Delta) \end{pmatrix} \in \mathbb{R}^{2k \times 2k},$$

$$G(\Delta) = \begin{pmatrix} 1 & 0 \\ -1/\Delta & 1/\Delta \end{pmatrix} \in \mathbb{R}^{2 \times 2},$$

(here $I_s$ is the $s \times s$ identity matrix) it is easy to verify that

(A.6) $$f(t,b) = L(\Delta)\psi(t,\tilde{b}_\Delta) + o(1),$$

where $\tilde{b}_\Delta = (b_1, b_1 + \Delta, \ldots, b_k, b_k + \Delta)^T$. For a fixed vector $T = (t_1, \ldots, t_m)^T \in \mathbb{R}^m$ with ordered components $t_1 < \cdots < t_m$ such that $t_i \in I$, $i = 1, \ldots, m$, define the matrices

$$F(T,b) = (f_i(t_j,b))_{i,j=1}^m, \qquad \psi(T,\tilde{b}) = (\psi_i(t_j,\tilde{b}))_{i,j=1}^m.$$

Then we obtain from (A.6)

(A.7)
$$\det F(T,b) = \lim_{\Delta \to 0} \frac{1}{\Delta^k} \psi(T, \tilde{b}_\Delta)$$
$$= \frac{\prod_{1 \leq i < j \leq m}(t_j - t_i) \prod_{1 \leq i < j \leq k}(b_i - b_j)^4}{\prod_{i=1}^k \prod_{j=1}^m (t_j - b_i)^2},$$

where the last identity follows from the fact that $\psi(T,\tilde{b})$ is a Cauchy–Vandermonde matrix, which implies

$$\det \psi(T,\tilde{b}) = \frac{\prod_{1 \leq i < j \leq m}(t_j - t_i) \prod_{1 \leq i < j \leq 2k}(\tilde{b}_i - \tilde{b}_j)}{\prod_{i=1}^{2k} \prod_{j=1}^{m}(t_j - \tilde{b}_i)}.$$

Now for any $b \in \mathcal{B}$, the right-hand side does not vanish and is of one sign independently of $T$. Consequently, $\{f_i(t,b)|i=1,\ldots,m\}$ is a Chebyshev system on the interval $I$. The assertion regarding the system $\{\bar{f}_i(t,x)|i=1,\ldots,m\}$ is proved similarly and, therefore, left to the reader. □

PROOF OF THEOREM 4.2. The second part of the theorem is a direct consequence of Lemma 4.1 and Theorem 3.2, while the first part of the proposition follows by Theorem A.2 in Imhof and Studden (2001). □

## APPENDIX B: AN AUXILIARY RESULT

Recall the notation in Sections 2 and 3, the definition of the regression functions in (2.1) and (3.3) and consider a design $\xi$ on the interval $I$ with at least $m$ support points. In this appendix we investigate the relation between the information matrices $M(\xi,b)$ and $\overline{M}(\xi,b)$ defined by (2.3) and (A.1), respectively, if condition (3.1) is satisfied, where the components of the vector $r = (r_1,\ldots,r_k)$ are different and ordered.

THEOREM B.1. *Assume that $\varphi \in C^{0,2k-1}$ and $\xi$ is an arbitrary design, such that the matrix $\overline{M}(\xi,b)$ is nonsingular. If assumption (3.1) is satisfied, it follows that for sufficiently small $\delta$ the matrix $M(\xi,b)$ is invertible and if $\delta \to 0$,*

$$(B.1) \quad M^{-1}(\xi,b) = \delta^{-4k+4} T(\delta) \begin{pmatrix} \overline{M}^{(1)}(\xi) & \overline{M}^{(2)}(\xi)F \\ F^T \overline{M}^{(2)^T}(\xi) & \gamma\gamma^T h + o(1) \end{pmatrix} T(\delta) + o(1),$$

*where the matrices $T(\delta) \in \mathbb{R}^{m \times m}$ and $\overline{M}^{(1)}(\xi) \in \mathbb{R}^{s \times s}, \overline{M}^{(2)}(\xi) \in \mathbb{R}^{s \times 2k}$ and $\overline{M}^{(3)}(\xi) \in \mathbb{R}^{2k \times 2k}$ are defined by*

$$T(\delta) = \operatorname{diag}\bigg( \underbrace{\delta^{2k-2},\ldots,\delta^{2k-2}}_{s}, \underbrace{\frac{1}{\delta},1,\frac{1}{\delta},1,\ldots,\frac{1}{\delta},1}_{2k} \bigg),$$

$$\begin{pmatrix} \overline{M}^{(1)} & \overline{M}^{(2)}(\xi) \\ \overline{M}^{(2)^T}(\xi) & \overline{M}^{(3)}(\xi) \end{pmatrix} = \overline{M}^{-1}(\xi,x),$$

*the vector $\gamma = (\gamma_1,\ldots,\gamma_{2k})^T$ and $h \in \mathbb{R}$ are given by $h = [(2k-1)!]^2 e_m^T \overline{M}^{-1}(\xi,x) e_m$,*

$$\gamma_{2i} = \prod_{j \neq i}(r_i - r_j)^{-2}, \qquad \gamma_{2i-1} = -\gamma_{2i} \sum_{j \neq i} \frac{2}{r_i - r_j}, \qquad i = 1,\ldots,k,$$

and the matrix $F \in \mathbb{R}^{2k \times 2k}$ is defined by

$$F = \begin{pmatrix} 0 & \cdots & 0 & \gamma_1/0! \\ \vdots & & & \\ 0 & \cdots & 0 & \gamma_{2k}/((2k-1)!) \end{pmatrix}.$$

PROOF. Define $\delta_i = r_i \delta$, $i = 1, \ldots, k$, $\psi(\delta) = (1, \delta, \ldots, \delta^{2k-1})^T$ and introduce the matrices

(B.2) $$L = (\ell_1, \ldots, \ell_{2k})^T \in \mathbb{R}^{2k \times 2k},$$

(B.3) $$U = \mathrm{diag}\left(1, \frac{1}{1!}, \frac{1}{2!}, \ldots, \frac{1}{(2k-1)!}\right) \in \mathbb{R}^{2k \times 2k},$$

where $\ell_{2i-1} = \psi(\delta_i)$, $\ell_{2i} = \psi'(\delta_i)$, $i = 1, \ldots, k$. For fixed $t \in I$, we use the Taylor expansions

$$\varphi(t, x + \delta) = \sum_{j=0}^{2k-1} \frac{\varphi^{(i)}(t, x)}{j!} \delta^j + o(\delta^{2k-1}),$$

$$\varphi'(t, x + \delta) = \sum_{j=1}^{2k-1} \frac{\varphi^{(i)}(t, x)}{(j-1)!} \delta^{j-1} + o(\delta^{2k-2}),$$

to obtain the representation

(B.4) $$f(t, b + \delta r) = \begin{pmatrix} I_s & 0 \\ 0 & LU \end{pmatrix} \bar{f}(t, x) + \begin{pmatrix} 0 \\ \tilde{f}(t) \end{pmatrix},$$

where $I_s \in \mathbb{R}^{s \times s}$ denotes the identity matrix and the vector $\tilde{f}$ is of order

(B.5) $$\tilde{f}(t) = (o(\delta^{2k-1}), o(\delta^{2k-2}), o(\delta^{2k-1}), \ldots, o(\delta^{2k-2}))^T.$$

It follows from pages 127–129 in Karlin and Studden (1966) that $\det L = \prod_{1 \leq i < j \leq k}(\delta_i - \delta_j)^4$ and consequently, $V = (v_1, \ldots, v_{2k}) := L^{-1}$ exists. The equality $LV = I_m$ implies the equations

$$v_{2i}^T \psi(\delta_j) = 0, \qquad v_{2i}^T \psi'(\delta_j) = 0, \qquad j \neq i,$$
$$v_{2i}^T \psi(\delta_i) = 0, \qquad v_{2i}^T \psi'(\delta_i) = 1,$$

which shows that $\delta_1, \ldots, \delta_{i-1}, \delta_{i+1}, \ldots, \delta_k$ are zeros of multiplicity two of the polynomial $v_{2i}^T \psi(\delta)$ and $\delta_i$ is a zero of multiplicity one. Because this polynomial has degree $2k - 1$, it follows that

(B.6) $$v_{2i}^T \psi(\delta) = (\delta - \delta_i) \prod_{j \neq i} \left(\frac{\delta - \delta_j}{\delta_j - \delta_i}\right)^2,$$

and a similar argument shows that

$$\text{(B.7)} \qquad v_{2i-1}^T \psi(\delta) = \frac{\delta - \alpha_i}{\delta_i - \alpha_i} \prod_{j \neq i} \left( \frac{\delta - \delta_j}{\delta_i - \delta_j} \right)^2,$$

where the constants $\alpha_1, \ldots, \alpha_k$ are given by

$$\text{(B.8)} \qquad \alpha_i = \delta_i + \left( \sum_{j \neq i} \frac{2}{\delta_i - \delta_j} \right)^{-1}, \qquad i = 1, \ldots, k.$$

From (B.4) and (B.5) we therefore obtain

$$f(t, b + \delta r) f^T(t, b + \delta r)$$
$$= \begin{pmatrix} I_s & 0 \\ 0 & LU \end{pmatrix} \bar{f}(t,x) \bar{f}^T(t,x) \begin{pmatrix} I_s & 0 \\ 0 & LU \end{pmatrix}^T + o(\delta^{2k-2}),$$

and integrating the right-hand side with respect to the design $\xi$ shows that

$$\text{(B.9)} \quad M(\xi, b + \delta r) = \begin{pmatrix} I_s & 0 \\ 0 & LU \end{pmatrix} \overline{M}(\xi, x) \begin{pmatrix} I_s & 0 \\ 0 & LU \end{pmatrix}^T + o(\delta^{2k-2}).$$

Now define $H_1(\delta) = \mathrm{diag}(\delta^{2k-1}, \delta^{2k-2}, \delta^{2k-1}, \ldots, \delta^{2k-1}, \delta^{2k-2}) \in \mathbb{R}^{2k \times 2k}$ and

$$H(\delta) = \begin{pmatrix} I_s & 0 \\ 0 & H_1(\delta) \end{pmatrix} \in \mathbb{R}^{m \times m}.$$

Then we obtain from (B.6) and (B.7) that $H_1(\delta)(L^{-1})^T = (0|\gamma) + o(1)$, where $\gamma = (\gamma_1, \ldots, \gamma_{2k})^T$ is defined by formula (B.2) and $0 \in \mathbb{R}^{2k \times 2k-1}$ denotes the matrix with all entries equal to zero. This implies that the inverse of the matrix $M(\xi, b + \delta r)$ is given by

$$M^{-1}(\xi, b + \delta r) = H^{-1}(\delta) \left\{ \begin{pmatrix} I & 0 \\ 0 & F \end{pmatrix} \overline{M}^{-1}(\xi, x) \begin{pmatrix} I & 0 \\ 0 & F^T \end{pmatrix} + o(1) \right\} H^{-1}(\delta)$$
$$= \delta^{-4k+4} T(\delta) \left\{ \begin{pmatrix} \overline{M}^{(1)}(\xi) & \overline{M}^{(2)}(\xi) F^T \\ F \overline{M}^{(2)^T}(\xi) & F \overline{M}^{(3)}(\xi) F^T \end{pmatrix} + o(1) \right\} T(\delta),$$

where the matrix $F$ is defined by $F = (0|\gamma) U^{-1} \in \mathbb{R}^{2k \times 2k}$. The assertion now follows by a straightforward calculation which shows that $F \overline{M}^{(3)}(\xi) F^T = h \gamma \gamma^T$.

$\square$

**Acknowledgments.** The authors are grateful to two unknown referees and the Associate Editor for their constructive comments on an earlier version of this paper, which led to a substantial improvement in the representation. The authors would also like to thank Isolde Gottschlich, who typed numerous versions of this paper with considerable technical expertise.

H. Dette
Fakultät für Mathematik
Ruhr-Universität Bochum
44780 Bochum
Germany
e-mail: holger.dette@ruhr-uni-bochum.de

V. B. Melas
A. Pepelyshev
St. Petersburg State University
Department of Mathematics
St. Petersburg
Russia
e-mail: v.melas@pobox.spbu.ru
e-mail: andrey@ap7236.spb.edu